\documentclass{amsart}
\usepackage{amssymb}
\usepackage{mathrsfs}
\usepackage{pifont}
\usepackage{amsfonts}
\usepackage{stmaryrd}
\usepackage{latexsym,amssymb,amsmath}

\newtheorem{theorem}{Theorem}[section]
\newtheorem{lemma}[theorem]{Lemma}
\newtheorem{corollary}[theorem]{Corollary}

\theoremstyle{definition}

\theoremstyle{remark}
\newtheorem{remark}[theorem]{Remark}

\numberwithin{equation}{section}

\begin{document}

\title{Abelian Ideals and Cohomology of Symplectic Type}

\author{Li Luo}
\address{Institute of mathematics, Academy of
Mathematics and System Sciences, Chinese Academy of Sciences,
Beijing 100080, China}
\email{luoli@amss.ac.cn}



\date{January 23, 2008.}


\keywords{Abelian ideal, cohomology, Symplectic Lie algebra, Weyl
group.}

\begin{abstract}
For symplectic Lie algebras $\mathfrak{sp}(2n,\mathbb{C})$, denote
by $\mathfrak{b}$ and $\mathfrak{n}$ its Borel subalgebra and
maximal nilpotent subalgebra, respectively. We construct a
relationship between the abelian ideals of $\mathfrak{b}$ and the
cohomology of $\mathfrak{n}$ with trivial coefficients. By this
relationship, we can enumerate the number of abelian ideals of
$\mathfrak{b}$ with certain dimension via the
Poincar$\acute{\mbox{e}}$ polynomials of Weyl groups of type
$A_{n-1}$ and $C_n$.
\end{abstract}

\maketitle


.

\section{Introduction}
Let $\mathfrak{g}$ be a finite dimension simple Lie algebra. Fixing
a Cartan subalgebra $\mathfrak{h}\subset\mathfrak{g}$, let
$\mathfrak{b}=\mathfrak{n}\oplus\mathfrak{h}$ be a Borel subalgebra
, where $\mathfrak{n}=[\mathfrak{b},\mathfrak{b}]$ is the maximal
nilpotent subalgebra of $\mathfrak{g}$. $\mathfrak{n}$ is also the
nilradical of $\mathfrak{b}$ under the Killing form. Needless to
say, the structure of the subalgebras $\mathfrak{b}$ and
$\mathfrak{n}$ is important to study $\mathfrak{g}$ itself.

The abelian ideals of the Boral subalgebra $\mathfrak{b}$, whose
origins can be traced back to work by Schur \cite{Si} (1905), have
recently enjoyed considerable attention. Kostant \cite{K2} (1998)
mentioned Peterson's $2^r$-theorem. That is, the number of abelian
ideals in $\mathfrak{b}$ is exactly $2^r$ where $r=\dim\mathfrak{h}$
is the rank of $\mathfrak{g}$. Moreover, Kostant found a relation
between abelian ideals of a Borel subalgebra and the discrete series
representations of the Lie group. Spherical orbits were described by
Panyushev and R\"{o}hrle \cite{PR} (2001) in terms of abelian
ideals. Furthermore, Panyushev \cite{P} (2003) discovered a
correspondence of maximal abelian ideals of a Borel subalgebra to
long positive roots. Suter \cite{S} (2004) determined the maximal
dimension among abelian subalgebras of a finite-dimensional simple
Lie algebra purely in terms of certain invariants and gave a uniform
explanation for Panyushev's result. Kostant \cite{K3} (2004) showed
that the powers of the Euler product and abelian ideals of a Borel
subalgebra are intimately related. Cellini and Papi \cite{CP} (2004)
had a detailed study of certain remarkable posets which form a
natural partition of all abelian ideals of a Borel subalgebra.

It is remarkable that the affine Weyl group $\widehat{\mathcal{W}}$
associated with $\mathfrak{g}$ is used in the proof of
$2^r$-theorem.

On the other hand, Bott \cite{B} gave a celebrated theorem which
shows that the Betti numbers $b_i$ of $\mathfrak{n}$ (i.e. the
dimension of $H^i(\mathfrak{n})$) can be expressed by the Weyl group
$\mathcal{W}$ associated with $\mathfrak{g}$. Later Kostant
\cite{K1} generalized his result to the nilradical of any parabolic
subalgebras $\mathfrak{p}$ of $\mathfrak{g}$. ($\mathfrak{n}$ is the
nilradical of $\mathfrak{b}$)

It seems that the Weyl group $\mathcal{W}$ (and its affine group
$\widehat{\mathcal{W}}$) can be a bridge to connect the cohomology
of $\mathfrak{n}$ to the abelian ideals of $\mathfrak{b}$. But no
one has given an explicit relationship between these two objects so
far. In this paper, we shall construct this relationship in the case
of $\mathfrak{g}=\mathfrak{sp}(2n,\mathbb{C})$.

In the following text, we let
$\mathfrak{g}=\mathfrak{sp}(2n,\mathbb{C})$ be the symplectic Lie
algebra. Our main theorem of this paper is as follows.

\vspace{0.2cm} \noindent {\bf Main Theorem} {\em
$H(\mathfrak{n})=\bigoplus_{(\sigma,I)\in S_n\times\mathcal
{I}}\mathbb{C}[L(\sigma,I)]$ where $S_n$ is the $n$-th symmetric
group and $\mathcal {I}$ is the set of all abelian ideals of
$\mathfrak{b}$. (The definition of $L: S_n\times \mathcal
{I}\rightarrow \wedge \mathfrak{n}^*$ is in (2.16). $[L(\sigma,I)]$
is the cohomology class defined by the harmonic cocycle
$L(\sigma,I)$.)}

\vspace{0.2cm} An interesting application of this theorem is to
compute the number of abelian ideals of $\mathfrak{b}$ with certain
dimension via the Poincar$\acute{\mbox{e}}$ polynomials of Weyl
groups of type $A$ and $C$. That is

\vspace{0.2cm} \noindent{\bf Corollary} {\em The number of abelian
ideals of $\mathfrak{b}$ with dimension $i$ is equal to the
coefficient of $t^i$ in $\prod_{i=1}^n(1+t^i)$.}

\section{Proof of the Main Result}
Let $\mathfrak{g}=\mathfrak{sp}(2n,\mathbb{C})$ be the symplectic
Lie algebra, $\mathcal{W}\simeq (\mathbb{Z}/2\mathbb{Z})^n\rtimes
S_n$ its Weyl group, $\Phi=\{\pm2\epsilon_i,
\pm(\epsilon_i\pm\epsilon_j)\mid 1\leq i\neq j\leq n\}$ its root
system with positive roots
$\Phi_+=\{2\epsilon_i,\epsilon_i\pm\epsilon_j\mid i<j\}$ and simple
roots $\pi=\{\epsilon_1-\epsilon_2, \epsilon_2-\epsilon_3,\ldots,
\epsilon_{n-1}-\epsilon_{n}, 2\epsilon_{n}\}$. Fixing a Cartan
subalgebra $\mathfrak{h}$, let
$\mathfrak{b}=\mathfrak{h}\bigoplus(\bigoplus_{\alpha\in
\Phi^+}\mathfrak{g}_\alpha)$ and
$\mathfrak{n}=[\mathfrak{b},\mathfrak{b}]=\bigoplus_{\alpha\in
\Phi^+}\mathfrak{g}_\alpha$ be the associated Borel subalgebra and
its nilradical, respectively.

Denote by $\mathcal {I}$ the set of all abelian ideals of
$\mathfrak{b}$. As mentioned in \cite{S}, there is a bijection as
follows.
\begin{eqnarray}
\Upsilon:=\{\Psi\subset \Phi_+\mid \Psi\dotplus \Phi_+\subset \Psi,
\Psi\dotplus \Psi=\emptyset\}&\leftrightarrow& \mathcal
{I}\nonumber\\\Psi&\mapsto&
I_\Psi:=\bigoplus_{\alpha\in\Psi}\mathfrak{g}_\alpha,
\end{eqnarray}
where $\Psi\dotplus \Phi_+:=(\Psi+\Phi_+)\cap\Phi_+$ and
$\Psi\dotplus \Psi:=(\Psi+\Psi)\cap\Phi_+$.

Denote $\Phi_+^0:=\{\epsilon_i-\epsilon_j\mid 1\leq i<j\leq n\}$.
\begin{lemma}
$\Psi\cap\Phi_+^0=\emptyset$ for any $\Psi\in\Upsilon$.
\end{lemma}
\begin{proof}
Suppose that there is an $\epsilon_i-\epsilon_j\in \Psi$. Since
$2\epsilon_j\in\Phi_+$ and
$(\epsilon_i-\epsilon_j)+2\epsilon_j=\epsilon_i+\epsilon_j\in\Phi_+$,
we have $\epsilon_i+\epsilon_j\in \Psi$. But
$(\epsilon_i-\epsilon_j)+(\epsilon_i+\epsilon_j)=2\epsilon_i\in\Phi_+$.
Contradiction.
\end{proof}

Define a partial ordering in
$$\Phi_+^1:=\Phi_+\setminus\Phi_+^0=\{\epsilon_i+\epsilon_j\mid 1\leq i\leq j\leq n\}$$ by
\begin{equation}
\epsilon_{i_1}+\epsilon_{j_1}\prec\epsilon_{i_2}+\epsilon_{j_2}
\Leftrightarrow i_1\geq i_2, j_1\geq j_2,
\end{equation}
where $i_1\leq j_1, i_2\leq j_2$.

A subset $\Psi\subset\Phi_+^1$ is called an increasing subset if for
any $x,y\in\Phi_+^1$, the conditions $x\in \Psi$ and $x\prec y$
imply $y\in \Psi$.

The following lemma is obvious.
\begin{lemma}
$\Upsilon = \mbox{the set of all increasing subsets of }\Phi_+^1$.
\hfill$\Box$
\end{lemma}

We shall show that these increasing subsets also appear on the
cohomology of $\mathfrak{n}$ with trivial coefficients.

Denote by $e_\alpha$ the unique element, up to nonzero scalar
multiples, of $\mathfrak{g}_\alpha$ for each $\alpha\in \Phi_+$.
Hence $\{e_\alpha\mid \alpha\in \Phi_+\}$ is a basis of
$\mathfrak{n}$. Define a linear function $f_\alpha\in
\mathfrak{n}^*$ by $f_\alpha(e_\beta)=\delta_{\alpha,\beta}$. Then
$\{f_\alpha\mid \alpha\in \Phi_+\}$ is a basis of $\mathfrak{n}^*$.

Recall the following theorem.
\begin{theorem}
({\bf Bott-Kostant.} c.f. \cite{B,K1})
\begin{equation}
H(\mathfrak{n})=\bigoplus_{w\in\mathcal{W}}\mathbb{C}[\wedge_{\alpha\in
\Phi_w}f_\alpha],\end{equation} where $\Phi_w=w(-\Phi_+)\cap\Phi_+$,
and $[\wedge_{\alpha\in \Phi_w}f_\alpha]$  is the cohomology class
defined by the (harmonic) cocycle $\wedge_{\alpha\in
\Phi_w}f_\alpha$. \hfill$\Box$
\end{theorem}

We see in Theorem 2.3 that $\Phi_w$ plays an important role in
cohomology. So we shall get some more information about $\Phi_w$.

Recall that for type $C_n$, the Weyl group $\mathcal{W}$ is
isomorphic to $(\mathbb{Z}/2\mathbb{Z})^n\rtimes S_n$. Precisely,
$\mathcal{W}$ can be realized by the composite of all permutations
of $\{1,2,\ldots,n\}$ and $\langle r_i\mid i=1,2,\ldots,n\rangle$
where $r_i(j)=-\delta_{i,j}i, (j=1,2,\ldots,n)$. In this paper, we
always use $(i_1,i_2,\ldots,i_n)$ to denote the permutation
$j\mapsto i_j (j=1,2,\ldots,n)$.

Each element in $\mathcal{W}$ can be expressed by the form
\begin{equation}
w=r_{j_1}r_{j_2}\cdots r_{j_k}(i_1,i_2,\ldots,i_n)(0\leq k\leq n),
\end{equation}
where
$(i_1,i_2,\ldots,i_n)(j_1)<(i_1,i_2,\ldots,i_n)(j_2)<\cdots<(i_1,i_2,\ldots,i_n)(j_k)$.
We call it the {\em standard form} of $w$.

\begin{lemma}
$\Phi_\sigma=\{\epsilon_i-\epsilon_j\mid 1\leq i<j\leq n,
\sigma^{-1}(i)>\sigma^{-1}(j)\}$ for any $\sigma\in
S_n\subset\mathcal{W}$. Hence if $\Phi_{\sigma_1}=\Phi_{\sigma_2}
(\sigma_1,\sigma_2\in S_n)$, then $\sigma_1=\sigma_2$.
\end{lemma}
\begin{proof}
It is obvious by a direct calculation.
\end{proof}

\begin{lemma}
For any $w\in\mathcal{W}$, there is a unique element $\eta_{w}\in
S_n$ such that $\Phi_{w}\cap\Phi_+^0=\Phi_{\eta_{w}}$. Precisely,
Write $w=r_{j_1}r_{j_2}\cdots r_{j_k}(i_1,i_2,\ldots,i_n)$ as the
standard form, then
\begin{equation}
\eta_{w}=(i_1,\ldots,\widehat{j_1},\ldots,\widehat{j_2},\ldots,\ldots,
\widehat{j_k},\ldots,i_n, j_k, j_{k-1},\ldots,j_1)
\end{equation} where $(i_1,\ldots,j_1,\ldots,j_2,\ldots,\ldots,
j_k,\ldots,i_n)=(i_1,i_2,\ldots,i_n)$ and the sign\quad
$\widehat{}$\quad indicates that the argument below it must be
omitted.
\end{lemma}

\begin{proof}
The uniqueness of $\eta_{w}$ comes from Lemma 2.4.

Denote
$\sigma=(i_1,\ldots,\widehat{j_1},\ldots,\widehat{j_2},\ldots,\ldots,
\widehat{j_k},\ldots,i_n, j_k, j_{k-1},\ldots,j_1)$.

We have to show $\Phi_{w}\cap\Phi_+^0=\Phi_\sigma$. In fact, for any
$\epsilon_i-\epsilon_j\in\Phi_{w}\cap\Phi_+^0$, it should be
$w^{-1}(\epsilon_i-\epsilon_j)=\epsilon_{w^{-1}(i)}-\epsilon_{w^{-1}(j)}\in-\Phi_+$.
What we need to do is to check $\sigma^{-1}(i)>\sigma^{-1}(j)$.
There are 4 cases as follows.

Case 1: $i,j\not\in\{j_1,j_2,\ldots,j_k\}$. Then
$w^{-1}(\epsilon_i-\epsilon_j)=\epsilon_{\sigma_0^{-1}(i)}-\epsilon_{\sigma_0^{-1}(j)}$,
where $\sigma_0=(i_1,i_2,\ldots,i_n)$. Hence
$\sigma_0^{-1}(i)>\sigma_0^{-1}(j)$. So
$\sigma^{-1}(i)>\sigma^{-1}(j)$ by the definition of $\sigma$.

Case 2: $i\in\{j_1,j_2,\ldots,j_k\},
j\not\in\{j_1,j_2,\ldots,j_k\}$. Then
$w^{-1}(\epsilon_i-\epsilon_j)=-\epsilon_{\sigma_0^{-1}(i)}-\epsilon_{\sigma_0^{-1}(j)}$
which is always in $-\Phi_+$.  We also do have
$\sigma^{-1}(i)>\sigma^{-1}(j)$ by the definition of $\sigma$.

Case 3: $i\not\in\{j_1,j_2,\ldots,j_k\},
j\in\{j_1,j_2,\ldots,j_k\}$. Then
$w^{-1}(\epsilon_i-\epsilon_j)=\epsilon_{\sigma_0^{-1}(i)}+\epsilon_{\sigma_0^{-1}(j)}$
which is always in $\Phi_+$. There should be no such $i,j$ with
$\sigma^{-1}(i)>\sigma^{-1}(j)$. It is the case by the definition of
$\sigma$.

Case 4: $i,j\in\{j_1,j_2,\ldots,j_k\}$. Then
$w^{-1}(\epsilon_i-\epsilon_j)=\epsilon_{\sigma_0^{-1}(j)}-\epsilon_{\sigma_0^{-1}(i)}$.
Hence $\sigma_0^{-1}(j)>\sigma_0^{-1}(i)$. So
$\sigma^{-1}(i)>\sigma^{-1}(j)$ by the definition of $\sigma$.
\end{proof}

By this lemma, we can define a map
\begin{equation}
\eta:\mathcal{W}\rightarrow S_n, \quad w\mapsto\eta_w.
\end{equation}

Denote \begin{equation} \sigma_l:=(n,n-1,\ldots,1)\in S_n,
\end{equation}
which is the longest element in $S_n$. It is clear that
$\sigma_l=\sigma_l^{-1}$ and
\begin{equation}
x\prec y\Leftrightarrow\sigma_l(y)\prec\sigma_l(x),\
(x,y\in\Phi_+^1).
\end{equation}

\begin{lemma}
$\sigma_l\eta_{w}^{-1}(\Phi_{w}\cap\Phi_+^1)\in \Upsilon.$
\end{lemma}

\begin{proof}
Write $w$ with the standard form $w=r_{j_1}r_{j_2}\cdots
r_{j_k}\sigma_0$, where $\sigma_0=(i_1,i_2,\ldots,i_n)$.

Take any $\epsilon_i+\epsilon_j\in\Phi_{w}\cap\Phi_+^1$. It should
be
$w^{-1}(\epsilon_i+\epsilon_j)=\epsilon_{w^{-1}(i)}+\epsilon_{w^{-1}(j)}\in-\Phi_+$.

There are 3 cases as follows.

Case 1: $i,j\not\in\{j_1,j_2,\ldots,j_k\}$. Then
$w^{-1}(\epsilon_i+\epsilon_j)=\epsilon_{\sigma_0^{-1}(i)}+\epsilon_{\sigma_0^{-1}(j)}$,
which can not be in $-\Phi_+$.

Case 2: Either $i$ or $j$, but not both, is in
$\{j_1,j_2,\ldots,j_k\}$. Without loss of generality, assume
$i\not\in\{j_1,j_2,\ldots,j_k\}$ and $j\in\{j_1,j_2,\ldots,j_k\}$.
Then
$w^{-1}(\epsilon_i+\epsilon_j)=\epsilon_{\sigma_0^{-1}(i)}-\epsilon_{\sigma_0^{-1}(j)}$,
which is in $-\Phi_+$ if and only if
$\sigma_0^{-1}(i)>\sigma_0^{-1}(j)$.

Case 3: $i,j\in\{j_1,j_2,\ldots,j_k\}$. Then
$w^{-1}(\epsilon_i+\epsilon_j)=-\epsilon_{\sigma_0^{-1}(i)}-\epsilon_{\sigma_0^{-1}(j)}$,
which is always in $-\Phi_+$.

All these three cases imply that for any $i$ and $j$ with
$\sigma_0(i)<\sigma_0(j)$, $\epsilon_i+\epsilon_j\in
\Phi_{w}\cap\Phi_+^1$ if and only if $i\in \{j_1,j_2,\ldots,j_k\}$.
So we have that for any $t\in\{1,2,\ldots,k\}$, if
$\epsilon_{j_t}+\epsilon_j\in\Phi_{w}\cap\Phi_+^1$, then (1)
$\epsilon_{j_t}+\epsilon_{j'}\in\Phi_{w}\cap\Phi_+^1$ for all $j'$
with $\sigma_0(j')>\sigma_0(j)$; and (2)
$\epsilon_{j_{t'}}+\epsilon_j\in\Phi_{w}\cap\Phi_+^1$ for all
$t'<t$. Hence by (2.5), we have that if
$\eta_{w}^{-1}(x)\in\Phi_{w}\cap\Phi_+^1$, then
$\eta_{w}^{-1}(y)\in\Phi_{w}\cap\Phi_+^1$ for any $y\in\Phi_+^1$
with $\eta_{w}^{-1}(y)\prec\eta_{w}^{-1}(x)$. So
$\sigma_l\eta_{w}^{-1}(\Phi_{w}\cap\Phi_+^1)\in \Upsilon$ by (2.8).
\end{proof}

\begin{remark}
In fact, the  proof of Lemma 2.6 also implies the explicit
expression of $\sigma_l\eta_{w}^{-1}(\Phi_{w}\cap\Phi_+^1)$. That is
for any $w=r_{j_1}r_{j_2}\cdots r_{j_k}\sigma_0\in \mathcal{W}\
(\sigma_0\in S_n)$,
\begin{equation}
\sigma_l\eta_{w}^{-1}(\Phi_{w}\cap\Phi_+^1)=\{\epsilon_i+\epsilon_j\mid
1\leq i\leq k, i\leq j\leq n+1-\sigma_0^{-1}(i)\}.
\end{equation}\hfill$\Box$
\end{remark}

Therefore we can define a map
\begin{equation}
\xi:\mathcal{W}\rightarrow\Upsilon,\quad
w\mapsto\sigma_l\eta_{w}^{-1}(\Phi_{w}\cap\Phi_+^1).
\end{equation}

\begin{lemma}
There exists a one to one correspondence as follows.
\begin{equation}
\mathcal{W}\leftrightarrow S_n\times\Upsilon,\quad w\mapsto(\eta(w),
\xi(w)).
\end{equation}
\end{lemma}

\begin{proof}
(2.5) and (2.9) imply that $w\mapsto(\eta(w), \xi(w))$ is injective.

Take any $(\sigma,\Psi)\in S_n\times\Upsilon$ and consider
$\sigma\sigma_l\Psi\subset \Phi_+^1$. We are going to determine the
$w\in\mathcal{W}$ such that $(\eta(w), \xi(w))=(\sigma,\Psi)$.

Assume $\{j\mid
2\epsilon_j\in\sigma\sigma_l\Psi\}=\{j_i,j_2,\ldots,j_k\}$ with
$\sigma^{-1}(j_1)>\sigma^{-1}(j_2)>\cdots>\sigma^{-1}(j_k)$. In
other words, $j_t=\sigma\sigma_l(t)\ (t=1,2,\ldots,k)$. Hence
$\sigma$ must be with the form
$\sigma=(i_1,i_2,\ldots,i_{n-k},j_k,j_{k-1},\ldots,j_{1})$ where
$\{i_1,i_2,\ldots,i_{n-k}\}$ is a permutation of
$\{1,2,\ldots,n\}\backslash\{j_i,j_2,\ldots,j_k\}$. Assume
$m_t=\max\{m\mid \epsilon_t+\epsilon_m\in\Psi\} (t=1,2,\ldots,k).$
It is clear that $m_1\geq m_2\geq\cdots\geq m_k$.

Take
\begin{eqnarray*}
w=(i_1,\ldots,i_{n-m_1},j_1,i_{n-m_1+1},\ldots,i_{n-m_k},j_k,i_{n-m_k+1},\ldots,i_{n-k})
\end{eqnarray*}
where
\begin{eqnarray*}
&&(\ldots,i_{n-m_t},j_t,i_{n-m_t+1},\ldots,i_{n-m_s},j_s,i_{n-m_s+1},\ldots)\\
&&:=(\ldots,i_{n-m_t},j_t,j_{t+1},\ldots,j_s,i_{n-m_s+1},\ldots)
\end{eqnarray*}
if $m_{t-1}<m_t=m_{t+1}=\cdots=m_s<m_{s+1}$. Then we can check
easily that $(\eta(w), \xi(w))=(\sigma,\Psi)$. So $w\mapsto(\eta(w),
\xi(w))$ is also surjective.
\end{proof}

\begin{remark}
By (2.1) and (2.11), we can obtain $|\mathcal
{I}|=|\Upsilon|=\frac{|(\mathbb{Z}/2\mathbb{Z})^n\rtimes
S_n|}{|S_n|}=2^n$ immediately. This is Peterson's $2^r$ theorem for
type $C$. \hfill$\Box$
\end{remark}

Each $\sigma\in S_n\subset\mathcal{W}$ induces a linear transform on
$\mathfrak{n}^*$ by
\begin{equation}
\sigma(f_{\epsilon_i\pm\epsilon_j})=f_{\epsilon_{\sigma(i)}\pm\epsilon_{\sigma(j)}},
\end{equation}
Moreover, this map can extend to $\bigwedge \mathfrak{n}^*$ by
\begin{equation}
\sigma(f_1\wedge f_2\wedge\cdots\wedge f_k)=\sigma(f_1)\wedge
\sigma(f_2)\wedge\cdots\wedge\sigma(f_k).
\end{equation}

By Lemma 2.4, (2.7), (2.10), (2.12) and $\sigma_l^{-1}=\sigma_l$, we
have
\begin{equation} \mathbb{C}\wedge_{\alpha\in
\Phi_w}f_\alpha=\mathbb{C}(\wedge_{\alpha\in\Phi_{\eta_w}}f_\alpha)
\wedge\eta_w\sigma_l(\wedge_{\alpha\in\xi_{w}}f_\alpha).
\end{equation}

Combine (2.1) and (2.14), then we get
\begin{equation}
\mathbb{C}\wedge_{\alpha\in
\Phi_w}f_\alpha=\mathbb{C}(\wedge_{\alpha\in\Phi_{\eta_w}}f_\alpha)
\wedge\eta_w\sigma_l(\wedge^{\max}I_{\xi_{w}}^*),
\end{equation}
where $I_{\xi_{w}}^*\subset\mathfrak{n}^*$ is the set of all linear
functions on $I_{\xi_{w}}$ and $\wedge^{\max}I_{\xi_{w}}^*$ is the
unique element (up to nonzero scalar multiples) in $\wedge^{\dim
I_{\xi_{w}} }I_{\xi_{w}}^*\subset \mathfrak{n}^*$.

Define
\begin{equation}
L: S_n\times \mathcal {I}\rightarrow \wedge \mathfrak{n}^*,
\quad(\sigma, I)\mapsto (\wedge_{\alpha\in\Phi_\sigma} f_\alpha)
\wedge\sigma\sigma_l(\wedge^{\max}I^*).
\end{equation}

Hence we can obtain our main theorem by (2.1), (2.3), (2.11),(2.14)
and (2.16) as follows.

\begin{theorem}
\begin{equation}
H(\mathfrak{n})=\bigoplus_{(\sigma,I)\in S_n\times\mathcal
{I}}\mathbb{C}[L(\sigma,I)]
\end{equation} where $[L(\sigma,I)]$ is the cohomology
class defined by the (harmonic) cocycle $L(\sigma,I)$.\hfill$\Box$
\end{theorem}

The definition of $L$ also implies that
\begin{equation}
\deg[L(\sigma,I)]=|\Phi_\sigma|+\dim I.
\end{equation}
Therefore
\begin{equation}
\dim H^i(\mathfrak{n})=\sum_{j+k=i}(|S_n^{(j)}|+|\mathcal
{I}^{(k)}|),
\end{equation}
where
\begin{equation}
S_n^{(j)}:=\{\sigma\in S_n\mid |\Phi_\sigma|=j\}.
\end{equation}
and
\begin{equation}
\mathcal {I}^{(k)}:=\{I\in \mathcal {I}\mid \dim I=k\}.
\end{equation}

On the other hand, Theorem 2.3 implies
\begin{equation}
\dim H^i(\mathfrak{n})=|((\mathbb{Z}/2\mathbb{Z})^n\rtimes
S_n)^{(i)}|
\end{equation}
where
\begin{equation}
((\mathbb{Z}/2\mathbb{Z})^n\rtimes S_n)^{(i)}=\{w\in
(\mathbb{Z}/2\mathbb{Z})^n\rtimes S_n\mid |\Phi_w|=i\}.
\end{equation}

The generate function of $|((\mathbb{Z}/2\mathbb{Z})^n\rtimes
S_n)^{(i)}|$\ and  $|S_n^{(i)}|$ are just the so-called
Poincar$\acute{\mbox{e}}$ polynomial of Weyl groups of type $C_n$
and $A_{n-1}$, respectively. These two polynomials can be found in
\cite{H}. We list them below.
\begin{equation}
\sum_{i=0}^{\infty}|((\mathbb{Z}/2\mathbb{Z})^n\rtimes
S_n)^{(i)}|t^i=\frac{\prod_{i=1}^n(1-t^{2i})}{(1-t)^n};
\end{equation}
\begin{equation}
\sum_{i=0}^{\infty}|S_n^{(i)}|t^i=\frac{\prod_{i=1}^{n}(1-t^{i})}{(1-t)^{n}}
\end{equation}

Thanks to (2.19) and (2.22), we get
\begin{equation}
\sum_{i=0}^\infty|\mathcal
{I}^{(i)}|t^i=\frac{\quad\frac{\prod_{i=1}^n(1-t^{2i})}{(1-t)^n}\quad}
{\quad\frac{\prod_{i=1}^{n}(1-t^{i})}{(1-t)^{n}}\quad}=\prod_{i=1}^n(1+t^i).
\end{equation}
That is
\begin{corollary}
The number of abelian ideals of $\mathfrak{b}$ with dimension $i$ is
equal to the coefficient of $t^i$ in
$\prod_{i=1}^n(1+t^i)$.\hfill$\Box$
\end{corollary}

\bibliographystyle{amsplain}

\end{document}